\documentclass[10pt]{elsarticle}
\journal{Applied Mathematics Letters}
\usepackage{lineno,hyperref}
\usepackage{color}
\usepackage{amsfonts,enumerate,amsmath,amssymb}
\def\qed{\strut\hfill $\Box$}
\newtheorem{thm}{Theorem}[section]

\def\para#1{\vskip .4\baselineskip\noindent{\bf #1}}
\numberwithin{equation}{section}
\begin{document}
	\begin{frontmatter}

		\title{The central limit theorem for slow-fast systems with L\'evy noise}
		
		\author[mymainaddress]{Xiaoyu Yang}
		\ead{yangxiaoyu@yahoo.com}
		
		\author[mymainaddress,mythirdaddress]{Yong Xu}
		\ead{hsux3@nwpu.edu.cn}

		\author[mymainaddress]{Ruifang Wang}
		\ead{wrfjy@yahoo.com}

		\author[mymainaddress]{Zhe Jiao\corref{mycorrespondingauthor}}
		\cortext[mycorrespondingauthor]{Corresponding author}
		\ead{zjiao@nwpu.edu.cn}

		\address[mymainaddress]{Department of Mathematics and Statistics, Northwestern Polytechnical University, Xi'an, 710072, China}
		\address[mythirdaddress]{MIIT Key Laboratory of Dynamics and Control of Complex Systems, Northwestern Polytechnical University, Xi'an, 710072, China}

		\begin{abstract}
			We consider a slow-fast stochastic differential system with  L\'evy noise. We will employ the perturbed test function method to study the normal deviation of the slow-fast system. Our main result states that the deviation can be approximated by a Gaussian process and the central limit theorem is obtained for the system. 
		\end{abstract}		
	\end{frontmatter}

	\section{Introduction}\label{sec-1}
	
	The slow-fast systems described by stochastic differential equations can be used to model many natural phenomena, such as weather and climate prediction \cite{Palmer}, macromolecules \cite{Qian}, geophysical flows \cite{Majda2006} and so on. We are concerned with the following slow-fast system driven by  L\'evy noise
\begin{equation} \label{orginal1.1}
\left
    \{
        \begin{array}{ll}
            d{x^\varepsilon }\left( t \right) = f_1\left( {{x^\varepsilon \left( t \right)},{y^\varepsilon \left( t \right)}} \right)dt+ {\sigma _1}dW_1\left( t \right) + \int_{\left| z \right| < 1} {{k _1}z\tilde N_1\left( {dz,dt} \right)},\\ 
             d{y^\varepsilon }\left( t \right) = \frac{1}{\varepsilon } {f_2}\left( {{y^\varepsilon \left( t \right)}} \right)dt + \frac{{\sigma _2}}{{\sqrt\varepsilon}}dW_2\left( t \right) + \int_{\left|z\right|<1}{{k _2} z\tilde N_2^{\varepsilon}\left({dz,dt}\right)},\\
             x^{\varepsilon}\left(0\right)=x_0 \in \mathbb{R},\quad {y^\varepsilon}\left(0\right)=y_0 \in \mathbb{R},
        \end{array} 
\right.
\end{equation}
for $t\in [0, T]$, $T<\infty$.
Here, $f_1: \mathbb{R}\times\mathbb{R} \rightarrow \mathbb{R}$, $f_2: \mathbb{R} \rightarrow \mathbb{R}$ are both nonlinear functions. 
The independent Poisson random measures $\tilde{N}_1(\cdot, \cdot)$ and $\tilde{N}_2^\varepsilon(\cdot, \cdot)$, are given by $\tilde{N}_1(dz, dt)=N_1(dz, dt)-\nu_1(dz) dt$, and
$\tilde{N}_2^\varepsilon(dz, dt)=N_2(dz, dt)-\frac{1}{\varepsilon}\nu_2(dz) dt$ respectively,
where $N_i(dz, dt)$, $i=1,2$, is associated Poisson measure, and $\nu_i$, $i=1,2$, the L\'evy measure satisfying $\int_\mathbb{R}(1\wedge z^2)\nu_i(dz)< \infty$. 
$W_i(t)$, $i=1,2$, are independent Brownian motions, which are also independent of $N_{i}(\cdot, \cdot)$. 
$\sigma_{i}$ and $k_i$ are all constants for $i=1,2$. The scaling parameter $\varepsilon>0$ is used to describe the separation of time scale between the slow variable $x_t^\varepsilon$ and the fast variable $y_t^\varepsilon$. 

Under some suitable conditions, it deduces that for $\varepsilon$ small enough, the trajectory of the slow variable $x_t^\varepsilon$ is situated in a small neighborhood of the function $\bar{x}_t$ which is the solution of the associated averaged system. This has attracted many researchers' attention (cf. \cite{Givon2007, Xu2011} and the references therein). It is interesting to study the deviation between $x_t^\varepsilon$ and $\bar{x}_t$.			
In this paper, we will show that  the difference between $x_t^\varepsilon$ and $\bar{x}_t$ weakly converges to a Gaussian process as $\varepsilon$ goes to zero.  For a related work on the normal deviation of the slow-fast stochastic differential systems, we refer to \cite{ WangRobert2013, LiuZhao2021}, which just consider the cases without the effect of L\'evy noise.  {Different from the literature focusing on the  Brownian motion, this is a new result owing to that in our situation, the tightness of deviation can not be proved by Ascoli-Arzela theorem directly.}
	
The plan of this paper is as follows. In the next section, we establish notation, give some precise conditions for the slow-fast system (\ref{orginal1.1}) and review some preliminary results. In Section \ref{sec-3}, we state and prove our main results. Throughout this paper, $c$, $C$, $c_1$, $C_1$, $\cdots$ denote certain positive constants that may vary from line to line.

	\section{Preliminaries}\label{sec-2}
Throughout this paper the quadruple $(\Omega, \mathcal{F}, \{\mathcal{F}_t\}_{t \geq 0},\mathbb{P})$ is a given stochastic basis satisfying the usual hypotheses. 
$\mathbb{E}(\cdot)$ stands for expectation with respect to the probability measure $\mathbb{P}$.	 
Let $\mathbb{D}=D\left([0,T],\mathbb{R}\right)$ be the space of $\mathbb{R}$-valued, right-continuous functions with left limits and the usual Skorohod topology. 	
The space $C^m (\mathbb{R}, \mathbb{R})$ consists of those functions that are $m-$times continuously differentiable, and $C^m_0 ([0,T], \mathbb{R})$ denote the subset of $C^m (\mathbb{R}, \mathbb{R})$ consisting of functions with compact support.

		\begin{itemize}
		\item[\textbf{A1}.] (Lipschitz continuity) For  any $ (x_1, y_1) $,  $ (x_2, y_2)\in \mathbb{R} \times\mathbb{R}$, there exists a constant $C_1> 0$ such that
		\begin{eqnarray*}
		\begin{aligned}
			&|f_1( x_2,y _2) - f_1( x_1,y _1)|^2 +|f_2( x_2,y _2) - f_2( x_1,y _1)|^2 \leq C_1(|x_2-x_1|^2+|y_2-y_1|^2).
		\end{aligned}
		\end{eqnarray*}
		\item[\textbf{A2}.] (Growth) For  all $(x, y)  \in \mathbb{R} \times\mathbb{R}$,  there exists a constant $ C_2>  0$ such that
		\begin{eqnarray*}
			|f_1( x, y)|^{2}+|f_2( x, y)|^{2}\leq C_2 (1+|x|^2+|y|^2).
		\end{eqnarray*}	
		\item[\textbf{A3}.] (Monotonicity) For  all $y_1, y_2\in \mathbb{R}$,  there exists a constant $ C_3>  0$ such that
		\begin{eqnarray*}
			(y_2 -y_1, f_2( y_2)-f_2( y_1))\leq -C_3 (|y_2-y_1|^2).
		\end{eqnarray*}	
		\item[\textbf{A4}.] (Regularity) For all $(x, y)  \in \mathbb{R} \times\mathbb{R}$, $\frac{\partial}{\partial x }f_1(x, y)$, $\frac{\partial^2}{\partial x^2 }f_1(x, y)$, and $\frac{\partial^3}{\partial x^3 }f_1(x, y)$ are bounded preserving with respect to $x$ and $y$. 
	\end{itemize}
Under assumptions (A1) to (A3), it implies from \cite[Theorem III.2.3.2]{JS1987} that there exists a unique strong solution of the following equation
\begin{eqnarray}\label{fastvariable}
d{y}\left( t \right) = f_2\left( {{y}\left( t \right)} \right)dt + {\sigma _2}dW_2\left( t \right) + \int_{\left| z \right| < 1} {z{k _2}\tilde N_2\left( {dz,dt} \right)}, \quad {y}\left( 0 \right)={y_0}.
\end{eqnarray}
Moreover, there also exists exactly one invariant measure $\mu$ with respect to the transition semigroup of $y(t)$, which is exponential mixing.
From \cite[Chapter 6]{Applebaum2009L}, assumptions (A1), (A2) and (A3) shows that there exists a unique mild solution $(x^\varepsilon(t), y^\varepsilon(t)) \in \mathbb{D} \times\mathbb{D} $ of system (\ref{orginal1.1}), and the slow variable satisfies 
 {\begin{eqnarray}\label{average}
		\mathbb{E}\big( {\mathop {\sup }\limits_{0 \le t \le T} \left| {{x^\varepsilon }\left( t \right) - \bar{x}\left( t \right)} \right|^2} \big) \le (\textrm{const})(\varepsilon/\delta),
\end{eqnarray}}
with suitable $\delta$ corresponding to $\varepsilon$ such that $\varepsilon/\delta \to 0$ when $\varepsilon \to 0$, and $\bar{x}_t$ is the solution of the averaged equation
\begin{eqnarray}\label{averagedeqn}
d\bar x\left( t \right) = \bar f_1\left( {\bar x\left( t \right)} \right)dt + {\sigma _1}dW_1\left( t \right) + \int_{\left| z \right| < 1} {{k _1}z\tilde N_1\left( {dz,dt} \right)}, 
\end{eqnarray}
where ${\bar f }_1\left( {\cdot} \right)=\int_{\mathbb{R}} {{f}_1\left( {{\cdot},y } \right){\mu }} \left( {dy } \right)$. 

Let $\mathcal{F}^{\varepsilon}_{t}$ denote the minimal $\sigma$-algebra generated by $\{x^\varepsilon(t), s\leq t\}$, and let $\mathbb{E}_t^\varepsilon$ be conditioning expectation on $\mathcal{F}^{\varepsilon}_{t}$.  Let $\mathcal{M}$ denote the set of real-valued progressively measurable functions $M: \Omega\times \mathbb{R} \rightarrow \mathbb{R}$ that are nonzero only on a bounded interval, and its subset denoted by ${\mathcal{M^\varepsilon} }$ whose element $M(t)$ is right-continuous, and also satisfies  
$\mathop {\sup }\limits_t \mathbb{E}|M \left( t \right)| < \infty$ and $M(t)$ is $\mathcal{F}_t^\varepsilon$-measurable. Let $M(t)$ and $M^{\delta}(t)$ belong to ${\mathcal{M^\varepsilon} }$ for each $\delta>0$. We say that $M(t)=p-\lim_{\delta\rightarrow 0} M^{\delta}(t)$ if and only if 
 	\begin{eqnarray*}
		\mathop {\sup }\limits_{t,\delta } \mathbb{E}\left| {{M ^\delta }\left( t \right)} \right| < \infty ,\quad\textrm{and} \quad 
		\mathop {\lim }\limits_{\delta  \to 0} \mathbb{E}\left| {{M^\delta }\left( t \right) - M\left( t \right)} \right| = 0, \quad \textrm{for each $t$}.
	\end{eqnarray*}	
 Define $\mathbf{T}(s)$ as a linear operator on ${\mathcal{M^\varepsilon} }$, and $\mathbf{T}(s)M(t)=\mathbb{E}_t^\varepsilon M(t+s) =\mathbb{E}(M(t+s)|\mathcal{F}_t^\varepsilon)$. It is not hard to show that $\mathbf{T}(s)$ is a semigroup. Then we can define the $p-$infinitesimal generator by
 \begin{eqnarray}\label{generator}
 	\mathcal{A}^{\varepsilon}M(t)=p-\lim_{\delta  \to 0}  {\frac{{\mathbf{T}(\delta)M\left( {t } \right) - M \left( t \right)}}{\delta }}
 \end{eqnarray}
 if the limit exists and be in ${\mathcal{M^\varepsilon} }$ \cite{Kushner1984Approximation}.

	\section{The Main Result}\label{sec-3}
 Then we define the normalized difference
\begin{eqnarray*}
	\Delta^{\varepsilon}(t) =\frac{1}{{\sqrt \varepsilon  }}\left({{x^\varepsilon }\left( t \right) - \bar{x}\left( t \right)}\right),
\end{eqnarray*}
which satisfies
\begin{eqnarray} \label{normaleqn}
	d\Delta^{\varepsilon}(t)= \frac{1}{{\sqrt \varepsilon  }}[ f_1( x^\varepsilon(t), y^{\varepsilon}(t))- \bar f_1(\bar{x}(t)) ]dt, \quad \Delta^{\varepsilon}(0)=0.
\end{eqnarray}
Our main theorem focuses on studying the weak convergence of the normal deviation $\Delta^{\varepsilon}(t) $, and proving that as $\varepsilon\rightarrow 0$, the limit of $\Delta^{\varepsilon}(t) $ weakly converges to  some kind of Gaussian process governed by the so-called linearized equation
\begin{eqnarray*}\label{linearized}
du(t) = \frac{\partial}{\partial \bar{x}}[\bar{f}_1(\bar{x}(t))] u(t)dt + H(\bar{x}( t )) dW_1(t), \quad u(0)=0,
\end{eqnarray*}	
where
\[
	{{H}\left( x \right)} = \Big(\int_0^\infty 2\mathbb{E} \big[{\big( {{f_1}( {x, \tilde{y}(s)} ) - {{\bar f_1}}\left( x \right)} \big)\big( {{f_1}( {x,\tilde{y}(0)} ) - {{\bar f_1}}\left( x \right)} \big)} \big]ds\Big)^{\frac{1}{2}}.
\]
Here, $\tilde{y}(t)$ is the solution of (\ref{fastvariable}) subject to the initial data $\tilde{y}(0)$ with distribution $\mu$, which deduces that the distribution of $\tilde{y}(t)$ is equal to $\mu$ for any $t\geq 0$. Thus, we give the statement of the main theorem as follows.
	\begin{thm}\label{thm}
		Let assumptions {\rm (A1)} to {\rm (A4)} hold. Then the normal deviation $\Delta^{\varepsilon}(t)$ converges to $u \left( t \right)$ in the sense of distribution as $\varepsilon$ goes to zero.
	\end{thm}
\para{Proof:} The proof is divided into four steps.	

\textbf{Step 1.} From (\ref{normaleqn}), we have
\[
		d\Delta^{\varepsilon}(t) = \frac{1}{{\sqrt \varepsilon  }}\big[ f_1( x^\varepsilon(t), y^{\varepsilon}(t))-f_1( {{ x^\varepsilon }\left( t \right),{\tilde{y}}( \frac{t}{\varepsilon})} ) +f_1( {{ x^\varepsilon }\left( t \right),{\tilde{y}}( \frac{t}{\varepsilon})} ) - \bar f_1\left( {\bar{x}\left( t \right)} \right) \big]dt,\quad \Delta^{\varepsilon}(0)=0.
\]
Then we can split $\Delta^{\varepsilon}(t) $ into $\Delta_1^{\varepsilon}(t)$ and $\Delta_2^{\varepsilon}(t)$ satisfying the following equations respectively
\begin{eqnarray}\label{linearized1}
		d\Delta_1^{\varepsilon}(t) = \frac{1}{{\sqrt \varepsilon  }}\big[f_1( x^\varepsilon(t), y^{\varepsilon}(t))-f_1( {{ x^\varepsilon }\left( t \right),{\tilde{y}}( \frac{t}{\varepsilon})} )\big]dt,\quad \Delta_1^{\varepsilon}(0)=0,
\end{eqnarray}
and
\begin{eqnarray}\label{linearized2}	
		d\Delta_2^{\varepsilon}(t) = \frac{1}{{\sqrt \varepsilon  }}\big[ f_1( {{ x^\varepsilon }\left( t \right),{\tilde{y}}( \frac{t}{\varepsilon})} ) - \bar f_1\left( {\bar{x}\left( t \right)} \right) \big]dt,\quad \Delta_2^{\varepsilon}(0)=0.
\end{eqnarray}	
From (\ref{linearized1}), (A1) and using H\"older's inequality, it follows that
	\begin{eqnarray}\label{linearized3}
	\begin{aligned}
		\mathbb{E} { \left| {\Delta_1 ^\varepsilon }\left( t \right) \right|} &=
		\mathbb{E} \Big[\frac{1}{{\sqrt \varepsilon  }}\int_0^t |f_1( x^\varepsilon(t), y^{\varepsilon}(t))-f_1(  x^\varepsilon(t),\tilde{y}( \frac{t}{\varepsilon}) )| ds \Big]\\
		&\le  \frac{C_1}{{ \sqrt \varepsilon  }}\int_0^t \Big[\mathbb{E}|y^{\varepsilon}(t)-{\tilde{y}}( \frac{t}{\varepsilon})|^2\Big]^{\frac{1}{2}}ds  \\
		&\le c_1\sqrt{\varepsilon}.
	\end{aligned}
	\end{eqnarray}	 		
If we can obtain that $\Delta_2^{\varepsilon}(t)$ converges to $u$ weakly as $\varepsilon\rightarrow 0$, then combining the estimate (\ref{linearized3}) and the following steps will show that the proof of the main theorem is completed. Therefore, in the following two steps we mainly investigate the convergence of $\Delta_2^{\varepsilon}(t)$.
	
\textbf{Step 2.}	
Define the truncated component ${\Delta_1 ^{\varepsilon ,K}}\left( t \right)$  which is the solution of the following equation
	\begin{eqnarray}\label{truncateeqn}
	d{\Delta_2 ^{\varepsilon ,K}}\left( t \right) =\frac{1}{{\sqrt \varepsilon  }}{q^{K}}( {{\Delta_2 ^{\varepsilon ,K}}\left( t \right)} )\big[ f( {{ x^\varepsilon }\left( t \right),{\tilde{y}}( \frac{t}{\varepsilon})} ) - \bar f\left( {\bar{x}\left( t \right)} \right) \big]dt, \quad \Delta_2^{\varepsilon, K}(0)=0,
	\end{eqnarray}	 
	where $K$ is a positive constant and 
	\begin{eqnarray*}
	{q^{K}}\big( {v}\big) =
	\left\{ 
	\begin{array}{l}
	1, \quad  v \in \left\{ {v :| v | \le K} \right\},\\
	0, \quad \textrm{otherwise}.
	\end{array}
	\right.
	\end{eqnarray*}	
We will prove that the solution ${\Delta_2 ^{\varepsilon,K} }\left( t \right)$ of (\ref{truncateeqn}) is tight in $\mathbb{D}$. Due to the Theorem 4 in \cite[Page 48]{Kushner1984Approximation}, we just need to show the following conditions hold.  
\begin{itemize}
	\item ${\Delta_2 ^{\varepsilon,K} }\left( t \right)$ is uniform bounded with probability 1, that is,
		\begin{eqnarray}\label{tight1}
			\mathop {\lim }\limits_{K \to \infty } \mathop {\lim \sup }\limits_{\varepsilon  \to 0} P\big\{ {\mathop {\sup }\limits_{0 \le t \le T} | {{\Delta_2 ^{\varepsilon ,K}}\left( t \right)} | \ge K} \big\} = 0.
		\end{eqnarray}	 
	\item There exists ${M^{\varepsilon,K} }\left(  t  \right) \in \textrm{Dom}(\mathcal{A}^{\varepsilon, K})$ such that 
		$\left\{ {\mathcal{A}^{\varepsilon,K}(M^{\varepsilon,K}\left( t \right)): t \le T} \right\}$
		is uniformly integrable, and 
		\begin{eqnarray}\label{tight2}
		p \text{-}{ \mathop {\lim }\limits_{\varepsilon  \to 0} \big( {{M^{\varepsilon,K} }\left(  t \right) - g( {{\Delta_2^{\varepsilon,K} }\left(  t \right)} )} \big)} = 0,
		\end{eqnarray}
		for each $g  \in C_0^2\left( {{\mathbb{R}},\mathbb{R}} \right)$. 
\end{itemize} 	
Indeed, it follows from the truncation that (\ref{tight1}) holds.   
And given a function $g  \in C_0^2\left( {{\mathbb{R}},\mathbb{R}} \right)$, we define 
\[
	{M ^{\varepsilon,K} }\left( t \right) = g( {{\Delta_2 ^{\varepsilon ,K}}\left( t \right)} ) + N_1^{\varepsilon,K} ( {{\Delta_2 ^{\varepsilon ,K}}\left( t \right)} )
\]
with
	\begin{eqnarray}\label{tight3}
	N_1^{\varepsilon,K} ( {{\Delta_2 ^{\varepsilon ,K}}\left( t \right)} ){\rm{ = }}{q^{K}}( {{\Delta_2 ^{\varepsilon ,K}}\left( t \right)} )\int_t^T {\frac{1}{{\sqrt \varepsilon  }}{g^{\prime}}( {{\Delta_2 ^{\varepsilon ,K}}\left( t \right)} )\mathbb{E}_t^\varepsilon \big[ f_1( {{ \bar{x}}\left( t \right),{\tilde{y}}( \frac{s}{\varepsilon})} ) - \bar f_1\left( {\bar{x}\left( t \right)} \right) \big]ds}.
	\end{eqnarray}
Here, $g^{\prime}(\cdot)$ and $g^{\prime\prime}(\cdot)$ denote the first-order derivative and the second-order derivative of $g$, respectively. We estimate (\ref{tight3})
	\begin{eqnarray*}
	\begin{aligned}
	\mathop {\sup }\limits_{0 \le t \le T} \mathbb{E}\big| {N_1^{\varepsilon,K} ( {{\Delta_1 ^{\varepsilon ,K}}\left( t \right)} )} \big| 
	&= \mathop {\sup }\limits_{0 \le t \le T} \sqrt \varepsilon  \mathbb{E}\big| {q^{K}}( {{\Delta_2 ^{\varepsilon ,K}}\left( t \right)} ){\int_{{t \mathord{\left/{\vphantom {t \varepsilon }} \right.\kern-\nulldelimiterspace} \varepsilon }}^{{T \mathord{\left/{\vphantom {T \varepsilon }} \right.\kern-\nulldelimiterspace} \varepsilon }} {{g^{\prime}}( {{\Delta_2 ^{\varepsilon ,K}}\left( t \right)} )\mathbb{E}_t^\varepsilon \big[ f_1( {{ \bar{x} }\left( t \right),{\tilde{y}}(\tau)} ) - \bar f_1\left( {\bar{x}\left( t \right)} \right) \big]d\tau} } \big|\\
	&\le \mathop {\sup }\limits_{0 \le t \le T} \sqrt \varepsilon  \mathbb{E}\Big| {\int_{{t \mathord{\left/{\vphantom {t \varepsilon }} \right.\kern-\nulldelimiterspace} \varepsilon }}^{{T \mathord{\left/{\vphantom {T \varepsilon }} \right.\kern-\nulldelimiterspace} \varepsilon }} {{g^{\prime}}( {{\Delta_2 ^{\varepsilon ,K}}\left( t \right)} )( {{e^{ -  \frac{\tau-t}{\varepsilon} }}} )d\tau} } \Big|\\
	&\le c_2{\sqrt \varepsilon  },
	\end{aligned}
	\end{eqnarray*}	
which implies (\ref{tight2}). Let $\mathcal{A}^{\varepsilon,K}$ be $p-$ infinitesimal generator associated to the truncated problem. Then we have
	\begin{eqnarray*}
	\mathcal{A}^{\varepsilon,K}{M ^{\varepsilon,K} }\left( t \right)=I_1 +I_2+I_3+I_4+I_5,
	\end{eqnarray*}	
where
\begin{eqnarray*}
\begin{aligned}
	I_1&=\frac{{q^{K}}}{\sqrt{\varepsilon}  }g^{\prime}( {{\Delta_2 ^{\varepsilon ,K}}\left( t \right)} )\big[ f_1( {{ x^\varepsilon }\left( t \right),{\tilde{y}}( \frac{t}{\varepsilon})} ) - f_1( {\bar{x}\left( t \right)}, {\tilde{y}}( \frac{t}{\varepsilon}) )\big],\\
	I_2 &= \frac{{q^{K}}}{{\varepsilon  }}\int_t^T { {{g^{\prime\prime}}( {\Delta_2 ^{\varepsilon ,K}( t )} )\mathbb{E}_t^\varepsilon \big[ f_1( {{ \bar{x} }\left( t \right),{\tilde{y}}( \frac{s}{\varepsilon})} ) - \bar f_1\left( {\bar{x}\left( t \right)} \right) \big]ds} } \big[ f_1( {{ x^\varepsilon }\left( t \right),{\tilde{y}}( \frac{t}{\varepsilon})} ) - \bar f_1\left( {\bar{x}\left( t \right)} \right) \big],\\
	I_3 &=\frac{{q^{K}}}{{\sqrt \varepsilon  }}\int_t^T {{{g^{\prime}}( {\Delta_2 ^{\varepsilon ,K}( t )} ){\frac{\partial}{\partial x}\big\{ {\mathbb{E}_t^\varepsilon \big[ f_1( {{ \bar{x} }\left( t \right),{\tilde{y}}( \frac{t}{\varepsilon})} ) -\bar f_1\left( {\bar{x}\left( t \right)} \right) \big]} \big\}}ds} }\times { \bar f_1( {{\bar{x }}\left( t \right)} )},\\
	I_4 &=\frac{\sigma_1^2{q^{K}}}{{2\sqrt \varepsilon  }}\int_t^T {{{g^{\prime}}( {\Delta_2 ^{\varepsilon ,K}( t )} ){\frac{\partial^2}{\partial x^2}\big\{ {\mathbb{E}_t^\varepsilon \big[ f_1( {{ \bar{x} }\left( t \right),{\tilde{y}}( \frac{t}{\varepsilon})} ) -\bar f_1\left( {\bar{x}\left( t \right)} \right) \big]} \big\}}ds} },\\
	I_5 &=\frac{{q^{K}}}{{\sqrt \varepsilon  }}\int_t^T\int_{\left| {z} \right| < 1}{g^{\prime}}( {\Delta_2 ^{\varepsilon ,K}( t )} )\Big\{{\mathbb{E}_t^\varepsilon \big[ f_1\big(( {{ \bar{x} }\left( t \right)+k_1z),{\tilde{y}}( \frac{t}{\varepsilon})} \big) -\bar f_1\left( {\bar{x}\left( t \right)} \right) \big]} \\
	&\qquad -{\mathbb{E}_t^\varepsilon \big[ f_1( {{ \bar{x} }\left( t \right),{\tilde{y}}( \frac{t}{\varepsilon})} ) -\bar f_1\left( {\bar{x}\left( t \right)} \right) \big]}-k_{1}z{\frac{\partial}{\partial x}\big\{ {\mathbb{E}_t^\varepsilon \big[ f_1( {{ \bar{x} }\left( t \right),{\tilde{y}}( \frac{t}{\varepsilon})} ) -\bar f_1\left( {\bar{x}\left( t \right)} \right) \big]} \big\}}\Big\}\nu_1\left( {dz} \right)ds.
\end{aligned}
\end{eqnarray*}					 
It is easy to check that $I_i$, $i=1, 2, 3, 4$, is uniformly integrable. By using Taylor's expansion, we have for $0<\theta<1$,
\[
	I_5=\frac{q^{K}k _1^2}{2\sqrt{\varepsilon}}\int_{\left| {z} \right| < 1} \int_t^T {{g^{\prime}}( {\Delta_2 ^{\varepsilon ,K}( t )} ){\frac{\partial^2}{\partial x^2}\big\{ {\mathbb{E}_t^\varepsilon [ {{f_1}\big( {(\bar{x}\left( t \right)+\theta z), {\tilde{y}}( \frac{t}{\varepsilon})} \big) - {{\bar f}_i}\left( \bar{x}\left( t \right)+\theta z \right)} ]} \big\}} }z^2ds  \nu_1( {dz}),
\]	
which implies that $I_5$ is uniformly integrable,						
then we prove that $\left\{ {\mathcal{A}^{\varepsilon,K}(M^{\varepsilon,K}\left( t \right)): t \le T} \right\}$ is uniformly integrable.  		

\textbf{Step 3.}	
	Since ${\Delta_2 ^{\varepsilon ,K}}\left( t \right)$ is tight in $\mathbb{D}$, then there exists a subsequence which converges to some $u^{K}(t)\in\mathbb{D}$ as $\varepsilon$ goes to zero. 
From (\ref{truncateeqn}), if follows that the limit $u^{K}(t)$ satisfies
   	\begin{eqnarray}\label{convergence}
   	du^K \left( t \right) =   \frac{\partial}{\partial x}[\bar f_1\left( {\bar{x}\left( t \right)} \right)]u^K \left( t \right)dt + {H}\left( {\bar{x}\left( t \right)} \right) dW_1.
   	\end{eqnarray}
And the corresponding infinitesimal operator for (\ref{convergence}) is given by
   	\begin{eqnarray*}
   	\mathcal{A}^{K} g \left( {u^K \left( t \right)} \right){\rm{ = }}g^{\prime}( {u^K \left( t \right)} )\frac{\partial}{\partial x}[\bar f_1\left( {\bar{x}\left( t \right)} \right)]u^K \left( t \right)
   	{\rm{ + }}\frac{1}{2}g^{\prime\prime} ( {u^K ( t )} )[{ H}\left( {\bar{x}( t )} \right)]^2.
   	\end{eqnarray*}	
If we prove that  for each more smooth function in $C_0^4\left( {{\mathbb{R}},\mathbb{R}} \right)$, denoted still by $g(\cdot)$, there exists ${L^{\varepsilon,K} }\left(  t  \right) \in \textrm{Dom}(\mathcal{A}^{\varepsilon, K})$ such that
		\begin{eqnarray}\label{convergencecond}
		\begin{aligned}
		&\qquad p-{\mathop {\lim }\limits_{\varepsilon  \to 0} [ {{L^{\varepsilon,K} }(  t ) - g( {{\Delta_2^{\varepsilon,K} }( t )} )} ]} = 0,\\
		&p-{\mathop {\lim }\limits_{\varepsilon  \to 0} [ {{\mathcal{A} ^{\varepsilon,K} }{L^{\varepsilon,K} }(  t ) - { \mathcal{A}^K}g( {{\Delta_2^{\varepsilon,K} }( t )} )} ]} = 0,
		\end{aligned}
		\end{eqnarray}
it implies from the Theorem 2 in \cite[Chapter 3]{Kushner1984Approximation} that   ${\Delta_1^{\varepsilon,K} }\left(  \cdot  \right)$ weakly converges to $u^K\left(  \cdot  \right)$ in $\mathbb{D}$ as $\varepsilon$ goes to zero.

To prove  (\ref{convergencecond}), we define
\[
	{L^{\varepsilon,K} }\left( t \right) = g( {{\Delta_2 ^{\varepsilon ,K}}\left( t \right)} ) + N_1^{\varepsilon,K} ( {{\Delta_2 ^{\varepsilon ,K}}\left( t \right)} )+N_2^{\varepsilon,K} ( {{\Delta_2 ^{\varepsilon ,K}}\left( t \right)} )
\]
with
	\begin{eqnarray*}
	N_2^{\varepsilon,K} ( {{\Delta_2 ^{\varepsilon ,K}}\left( t \right)} ){\rm{ = }}{q^{K}}\int_t^T {{\Delta_2 ^{\varepsilon ,K}(t)}{g^{\prime}}( {{\Delta_2 ^{\varepsilon ,K}}\left( t \right)} )\mathbb{E}_t^\varepsilon \big[ \frac{\partial}{\partial x}f_1( {{ \bar{x}}\left( t \right),{\tilde{y}}( \frac{s}{\varepsilon})} ) -\frac{\partial}{\partial x} \bar f_1\left( {\bar{x}\left( t \right)} \right) \big]ds}.
	\end{eqnarray*}
Note that $\mathop {\sup }\limits_{0 \le t \le T} \mathbb{E}\big| {N_2^{\varepsilon,K} ( {{\Delta_1 ^{\varepsilon ,K}}\left( t \right)} )} \big|\leq c_3\sqrt \varepsilon$. Then with (\ref{tight2}),  the first estimate of (\ref{convergencecond}) is obtained. Similarly, we have
	\begin{eqnarray*}
	\mathcal{A}^{\varepsilon,K}{N_2^{\varepsilon,K} }\left( t \right)=J_1 +J_2+J_3+J_4+J_5,
	\end{eqnarray*}	
where
\begin{eqnarray*}
\begin{aligned}
	J_1&=-{q^{K}}{{\Delta_2 ^{\varepsilon ,K}}\left( t \right)}g^{\prime}( {{\Delta_2 ^{\varepsilon ,K}}\left( t \right)} ) \big[ \frac{\partial}{\partial x}f_1( {{ \bar{x}}\left( t \right),{\tilde{y}}( \frac{t}{\varepsilon})} ) -\frac{\partial}{\partial x} \bar f_1\left( {\bar{x}\left( t \right)} \right) \big],\\
	J_2 &={q^{K}}\int_t^T { {[{{\Delta_2 ^{\varepsilon ,K}}\left( t \right)}{g^{\prime}}( {\Delta_2 ^{\varepsilon ,K}( t )} )]^{\prime}\mathbb{E}_t^\varepsilon \big[ \frac{\partial}{\partial x}f_1( {{ \bar{x} }\left( t \right),{\tilde{y}}( \frac{s}{\varepsilon})} ) - \frac{\partial}{\partial x}\bar f_1\left( {\bar{x}\left( t \right)} \right) \big]ds} } \big[ f_1( {{ x^\varepsilon }\left( t \right),{\tilde{y}}( \frac{t}{\varepsilon})} ) - \bar f_1\left( {\bar{x}\left( t \right)} \right) \big],\\
	J_3 &={q^{K}}\int_t^T {{{{\Delta_2 ^{\varepsilon ,K}}\left( t \right)}g^{\prime}( {{\Delta_2 ^{\varepsilon ,K}}\left( t \right)} ) {\frac{\partial}{\partial x}\big\{ {\mathbb{E}_t^\varepsilon \big[ \frac{\partial}{\partial x}f_1( {{ \bar{x} }\left( t \right),{\tilde{y}}( \frac{t}{\varepsilon})} ) -\frac{\partial}{\partial x}\bar f_1\left( {\bar{x}\left( t \right)} \right) \big]} \big\}}ds} }\times { \bar f_1( {{\bar{x }}\left( t \right)} )},\\
	J_4 &=\frac{\sigma_1^2{q^{K}}}{{2  }}\int_t^T {{{{\Delta_2 ^{\varepsilon ,K}}\left( t \right)}g^{\prime}( {{\Delta_2 ^{\varepsilon ,K}}\left( t \right)} ) {\frac{\partial^2}{\partial x^2}\big\{ {\mathbb{E}_t^\varepsilon \big[ \frac{\partial}{\partial x}f_1( {{ \bar{x} }\left( t \right),{\tilde{y}}( \frac{t}{\varepsilon})} ) -\frac{\partial}{\partial x}\bar f_1\left( {\bar{x}\left( t \right)} \right) \big]} \big\}}ds} },\\
	J_5 &={q^{K}}\int_t^T\int_{\left| {z} \right| < 1}{{\Delta_2 ^{\varepsilon ,K}}\left( t \right)}g^{\prime}( {{\Delta_2 ^{\varepsilon ,K}}\left( t \right)} ) \Big\{{\mathbb{E}_t^\varepsilon \big[ \frac{\partial}{\partial x}f_1\big(( {{ \bar{x} }\left( t \right)+k_1z),{\tilde{y}}( \frac{t}{\varepsilon})} \big) -\frac{\partial}{\partial x}\bar f_1\left( {\bar{x}\left( t \right)} \right) \big]} \\
	&\qquad -{\mathbb{E}_t^\varepsilon \big[ \frac{\partial}{\partial x}f_1( {{ \bar{x} }\left( t \right),{\tilde{y}}( \frac{t}{\varepsilon})} ) -\frac{\partial}{\partial x}\bar f_1\left( {\bar{x}\left( t \right)} \right) \big]}-k_{1}z{\frac{\partial}{\partial x}\big\{ {\mathbb{E}_t^\varepsilon \big[\frac{\partial}{\partial x} f_1( {{ \bar{x} }\left( t \right),{\tilde{y}}( \frac{t}{\varepsilon})} ) -\frac{\partial}{\partial x}\bar f_1\left( {\bar{x}\left( t \right)} \right) \big]} \big\}}\Big\}\nu_1\left( {dz} \right)ds.
\end{aligned}
\end{eqnarray*}	
By the assumption (A3) and the exponential mixing of the invariant measure, we can also obtain that $\sum_{i=3}^{5}I_i$ and $\sum_{j=2}^{5}J_j$ both converge to $0$ as $\varepsilon\rightarrow 0$. Moreover, we also have 
\begin{eqnarray}
\begin{aligned} \label{plimit1}
I_1+J_1&={q^{K}}g^{\prime}( {{\Delta_2 ^{\varepsilon ,K}}\left( t \right)} )\frac{f_1( {{ x^\varepsilon }\left( t \right),{\tilde{y}}( \frac{t}{\varepsilon})} ) - f_1( {\bar{x}\left( t \right)}, {\tilde{y}}( \frac{t}{\varepsilon}) )}{{ x^\varepsilon }\left( t \right)-{\bar{x}\left( t \right)}}{{\Delta_2 ^{\varepsilon ,K}}\left( t \right)}+J_1\\
&={q^{K}}{{\Delta_2 ^{\varepsilon ,K}}\left( t \right)}g^{\prime}( {{\Delta_2 ^{\varepsilon ,K}}\left( t \right)} )\frac{\partial}{\partial x}f_1( {{ \bar{x}}\left( t \right),{\tilde{y}}( \frac{t}{\varepsilon})} )+o(\sqrt{\varepsilon})+J_1\\
&={q^{K}}{{\Delta_2 ^{\varepsilon ,K}}\left( t \right)}g^{\prime}( {{\Delta_2 ^{\varepsilon ,K}}\left( t \right)} ) \frac{\partial}{\partial x} \bar f_1\left( {\bar{x}\left( t \right)} \right) +o(\sqrt{\varepsilon}),
\end{aligned}
\end{eqnarray}	
and $I_2 = I_{21}+I_{22}$ where
\begin{eqnarray*}
\begin{aligned}
	I_{21}&=q^{K}{g^{\prime\prime}}( {\Delta_2 ^{\varepsilon ,K}( t )} )\int_{\frac{t}{\varepsilon}}^{\frac{T}{\varepsilon}}\mathbb{E}_t^\varepsilon \big[ f_1( {{ \bar{x} }\left( t \right),{\tilde{y}}( \tau)} ) - \bar f_1\left( {\bar{x}\left( t \right)} \right) \big]d\tau[f_1( {{ x^\varepsilon }\left( t \right),{\tilde{y}}( \frac{t}{\varepsilon})} )-f_1( {{ \bar{x} }\left( t \right),{\tilde{y}}( \frac{t}{\varepsilon})} )],\\
	I_{22}&=q^{K}{g^{\prime\prime}}( {\Delta_2 ^{\varepsilon ,K}( t )} )\int_{\frac{t}{\varepsilon}}^{\frac{T}{\varepsilon}}\mathbb{E}_t^\varepsilon \big[ f_1( {{ \bar{x} }\left( t \right),{\tilde{y}}( \tau)} ) - \bar f_1\left( {\bar{x}\left( t \right)} \right) \big]d\tau[f_1( {{ \bar{x} }\left( t \right),{\tilde{y}}( \frac{t}{\varepsilon})} ) - \bar f_1\left( {\bar{x}\left( t \right)} \right)].
\end{aligned}
\end{eqnarray*}	 
By the exponential mixing of the invariant measure, we know $\mathbb{E}[I_{21}]=o(\varepsilon)$. Moreover, we have
\begin{eqnarray}\label{plimit2}
	\mathbb{E}[I_{22}]=\frac{1}{2}q^{K}{g^{\prime\prime}}( {\Delta_2 ^{\varepsilon ,K}( t )} ) [H({\bar{x}\left( t \right)})]^2+o(\varepsilon). 		 	
\end{eqnarray}
Due to $\mathcal{A}^{\varepsilon,K}[L^{\varepsilon,K} \left( t \right)]=\sum_{i=1}^{5}(I_i+J_i)$, it deduce from (\ref{plimit1}) and (\ref{plimit2}) that the second estimate of (\ref{convergencecond}) holds true.	
 
\textbf{Step 4.}			
Taking $K \to \infty $, together with  the weak convergence of ${\Delta_2 ^{\varepsilon ,K}}\left(  \cdot  \right)$, ${\Delta_2 ^\varepsilon }\left(  t  \right)$ converges weakly to $u \left(  t  \right)$  in  $\mathbb{D}$ as $\varepsilon$ goes to zero.

This proof is completed. 
\qed

	\section*{Acknowledgments}
	This work was partly supported by the Key International (Regional) Cooperative Research
	Projects of the NSF of China (Grant 12120101002), the NSF of China
	(Grant 12072264, 11802236), the Fundamental Research Funds for the Central Universities, the
	Research Funds for Interdisciplinary Subject of Northwestern Polytechnical University, 
	the Shaanxi Provincial Key R\&D Program
	(Grants 2020KW-013, 2019TD-010).


\section*{References}
	
	\begin{thebibliography}{}
		
\bibitem{Palmer}T. N. Palmer, A nonlinear dynamical perspective on model error: A proposal for non-local stochastic-dynamic parametrization in weather and climate prediction models, Quarterly Journal of the Royal Meteorological Society, 127, 279-304 (2001).

\bibitem{Qian}H. Qian, Mesoscopic nonequilibrium thermodynamics of single macromolecules and dynamic entropy-energy compensation. Physical Review E, 65, 016102 (2001).

\bibitem{Majda2006}A. Majda, X. Wang, Nonlinear dynamics and statistical theories for basic geophysical flows. Cambridge University Press, 2006.

		



\bibitem{Givon2007}D. Givon, Strong convergence rate for two-time-scale jump-diffusion stochastic differential systems, Multiscale Modeling and Simulation, 6,  577-594 (207).
		
		
\bibitem{Xu2011}Y. Xu, J. Duan, W. Xu, An averaging principle for stochastic dynamical systems with L\'evy noise, Physica D, 240,  1395-1401 (2011).

		
\bibitem{WangRobert2013}W. Wang, A.J. Roberts, Slow manifold and averaging for slow–fast stochastic differential system, Journal of Mathematical Analysis and Applications, 398, 822-839 (2013).
		
		
		\bibitem{LiuZhao2021}J. Liu, M. Zhao, Normal deviation of synchronization of stochastic coupled systems, Discrete \& Continuous Dynamical Systems-B, Online (2021).
		

		
		
\bibitem{Applebaum2009L}D. Applebaum, L\'evy Processes and Stochastic Calculus, Cambridge University Press, Second Edition, 2009.		
		
		
		
\bibitem{JS1987}J. Jacod and A. N. Shiryaev. Limit Theorems for Stochastic Processes. Springer, Berlin, 1987.
		
		
		
		\bibitem{Kushner1984Approximation}H. Kushner, Approximation and weak convergence methods for random processes, with applications to stochastic systems theory, MIT Press, Cambridge, MA, 1984.
		


		
	\end{thebibliography}
\end{document}